\documentclass[12pt,a4paper,reqno]{amsart}
\usepackage{amsmath}
\usepackage{amsfonts}
\usepackage{amssymb}
\numberwithin{equation}{section}

\theoremstyle{plain}

\newtheorem{theorem}[subsection]{Theorem}
\newtheorem{proposition}[subsection]{Proposition}
\newtheorem{lemma}[subsection]{Lemma}
\newtheorem{claim}[subsection]{Claim}
\newtheorem{corollary}[subsection]{Corollary}
\newtheorem{conjecture}[subsection]{Conjecture}

\theoremstyle{definition}

\newtheorem{definition}[subsection]{Definition}

\theoremstyle{remark}

\newtheorem*{example}{Example}
\newtheorem*{examples}{Examples}
\newtheorem*{remark}{Remark}
\newtheorem*{remarks}{Remarks}

\renewcommand{\leq}{\leqslant}
\renewcommand{\geq}{\geqslant}

\newsavebox{\proofbox}
\savebox{\proofbox}{\begin{picture}(7,7)%
  \put(0,0){\framebox(7,7){}}\end{picture}}
\def\boxeq{\tag*{\usebox{\proofbox}}}

\def\proof{\noindent\textit{Proof. }}
\def\endproof{\hfill{\usebox{\proofbox}}\vspace{9pt}}

\def\N{\mathbb{N}}
\def\F{\mathbb{F}}
\def\Z{\mathbb{Z}}
\def\R{\mathbb{R}}

\def\P{\mathbb{P}}

\def\tF{{F_{\operatorname{SMD}}}}
\def\tG{{G_{\operatorname{SMD}}}}
\def\dF{{F_{\operatorname{D}}}}
\def\dG{{G_{\operatorname{D}}}}
\def\half{\textstyle\frac{1}{2}\displaystyle}

\newcommand\HS{\operatorname{HS}}

\begin{document}

\title[Freiman in finite fields]{Freiman's theorem in finite fields via extremal set theory}

\author{Ben Green}
\address{Centre for Mathematical Sciences, Wilberforce Road, Cambridge CB3 0WA, England.}
\email{b.j.green@dpmms.cam.ac.uk}

\author{Terence Tao}
\address{Department of Mathematics, UCLA.
}
\email{tao@math.ucla.edu}

\thanks{The first author is a Clay Research Fellow, and is pleased to acknowledge the support of the Clay Mathematics Institute.  The second author is supported by a grant from the MacArthur Foundation.}

\begin{abstract} 
Using various results from extremal set theory (interpreted in the language of additive combinatorics), we prove an asyptotically sharp version of Freiman's theorem in $\F_2^n$: if $A \subseteq \F_2^n$ is a set for which $|A + A| \leq K|A|$ then $A$ is contained in a subspace of size $2^{2K + O(\sqrt{K}\log K)}|A|$; except for the $O(\sqrt{K} \log K)$ error, this is best possible. If in addition we assume that $A$ is a downset, then we can also cover $A$ by $O(K^{46})$ translates of a coordinate subspace of size at most $|A|$, thereby verifying the so-called polynomial Freiman-Ruzsa conjecture in this case.  A common theme in the arguments is the use of compression techniques. These have long been familiar in extremal set theory, but have been used only rarely in the additive combinatorics literature.
\end{abstract}

\maketitle

\section{Introduction and statement of results}

If $A,B$ are any two subsets of an abelian group $G$ we write $A + B := \{a + b : a \in A, b \in B\}$ for the collection of pairwise sums of $A$ and $B$, and $|A|$ for the cardinality of $A$. Define the \emph{doubling constant} $\sigma(A)$ of a non-empty finite set $A$ by $\sigma[A] := \frac{|A+A|}{|A|}$.  A fundamental theorem of Freiman \cite{freiman} in additive combinatorics asserts that non-empty finite sets of integers with small doubling constant have a special structure, namely that they can be efficiently contained in a generalized arithmetic progression.  Similar results are known for sets in any abelian group \cite{gr-4}, in which the notion of a generalized arithmetic progression is replaced by the slightly more general notion of a \emph{coset progression}.  See \cite[Chapter 5]{tao-vu} for further discussion.

The situation becomes particularly simple in the \emph{finite field model}, in which the ambient group $G$ has a fixed finite torsion.  In this paper we shall concern ourselves exclusively with the two-torsion case 
$$G = \F_2^I := \{ (x_i)_{i \in I}: x_i \in \F_2 \hbox{ for all } i \in I \}$$
where $\F_2$ is the finite field with two elements and $I \subseteq \Z_+$ is a finite set of natural numbers; we abbreviate\footnote{Of course, every $\F_2^I$ is isomorphic to $\F_2^n$ for some $n$, but to avoid some tedious relabeling later on in the arguments it is convenient to work in the slightly this slightly more general setting when the label set $I$ does not need to be an initial segment.} $\F_2^{\{1,\ldots,n\}}$ as $\F_2^n$.  If $J \subseteq I$ then we embed $\F_2^J$ into $\F_2^I$ in the obvious manner (extending by zeroes), and so in particular $\F_2^m \subseteq \F_2^n$ whenever $m < n$.  We refer to $\F_2^J$ as a \emph{coordinate subspace} of $\F_2^I$.  The article \cite{green-finite-field-survey} provides an extended discussion of the r\^ole of finite field models in additive combinatorics.  As we shall soon see, additive combinatorics in this setting is also connected to more classical combinatorial problems in set systems.

In the finite field model, the analogue of a generalized arithmetic progression or coset progression is an affine subspace (i.e. a translate of a linear subspace of $\F_2^I$); note that these are the unique sets of doubling constant $1$.  We thus expect that sets of small doubling constant to be related to affine subspaces.  To quantify this relationship, let us define two quantities.

\begin{definition}[$F(K)$ and $G(K)$]\label{fg-def} Let $K \geq 1$.  
We define $F(K)$ to be the least constant such that for any finite $I \subseteq \Z_+$ and any non-empty $A \subseteq \F_2^I$ with doubling constant $\sigma[A] \leq K$, there exists an affine subspace $V \subseteq \F_2^I$ of cardinality $|V|\leq F(K) |A|$ which contains $A$.  

We define $G(K)$ to be the least constant such that for any finite $I \subseteq \Z_+$ and any non-empty $A \subseteq \F_2^I$ with doubling constant $\sigma[A] \leq K$, there exists a linear subspace $V \subseteq \F_2^I$ of cardinality $|V| \leq |A|$ such that $A$ is covered by at most $G(K)$ translates of $V$.
\end{definition}

Thus for instance one can easily verify that $F(1)=G(1)=1$, and that $F$ and $G$ are both non-decreasing in $K$.  

It is not immediately obvious that $F$ and $G$ are even finite.  However, in \cite{ruzsa-freiman} the first non-trivial upper bound on $F$ was shown, namely $F(K) \leq K^2 2^{K^4}$ for all $K \geq 1$. This was improved by Ruzsa (see \cite{desh}) to $F(K) \leq K 2^{\lfloor K \rfloor^3 - 1}$, by the first author and Ruzsa \cite{green-ruzsa} to $F(K) \leq K^2 2^{2K^2-2}$ and subsequently by Sanders \cite{sanders} to $F(K) \leq 2^{O(K^{3/2}\log (1+K))}$. Some additional results for $F(K)$ are known when $K$ is small: see \S \ref{smallk}. However, our primary interest here is with the large $K$ case (e.g. $K \geq 2$).

A lower bound for $F(K)$ is not hard to obtain, as the following example shows.

\begin{lemma}[Independent points example]\label{ind-pt-example}  For all $K \geq 1$ we have
\[ F(K) \geq 2^{2K - O( \log (1+K) )}.\]
\end{lemma}
\proof Let $r := \lfloor 2K - 1\rfloor$, and consider the set $A = \{0,e_1,\dots,e_r\} \subseteq \F_2^r$. We have $|A + A| = \frac{1}{2}(r^2 + r + 2)$ and so $\sigma[A] \leq \frac{1}{2}(r + 1) \leq K$. Furthermore $A$ is not contained in any subspace of dimension less than $r$. Thus $F(K)|A| \geq 2^r$, and whence $\log_2 F(K) \geq r - \log_2 |A| = 2K - O(\log (1+K))$ as desired.\endproof

The first main theorem in this paper, which we prove in \S \ref{freisec}, provides an upper bound to match this trivial lower bound.  
\begin{theorem}[Asymptotic for $F(K)$]\label{mainthm}
For all $K \geq 2$ we have 
\[ F(K) = 2^{2K + O(\sqrt{K} \log K)}.\]
\end{theorem}
Except for the issue of improving the bounds on the error term $O(\sqrt{K} \log K)$, this settles the question of determining $F(K)$.  Now we turn to the quantity $G(K)$.  This quantity is related to $F(K)$ by the easily verified inequalities
$$ 1+\log_2 F(K) \leq G(K) \leq 2 F(K).$$
In particular, Theorem \ref{mainthm} gives $G(K) \leq 2^{2K + O(\sqrt{K} \log K)}$ for $K \geq 2$.  However, one can do better than this.  The following conjecture was attributed to Marton in \cite{ruzsa-freiman} and is sometimes known as the \emph{polynomial Freiman-Ruzsa conjecture} (PFR).
\begin{conjecture}[Polynomial Freiman-Ruzsa conjecture]\label{pfr} For any $K \geq 1$ we have the inequality $G(K) \ll K^{O(1)}$.  In other words, any non-empty set $A \subseteq \F_2^I$ is contained in the union of $O(K^{O(1)})$ cosets of a subspace of cardinality at most $|A|$.
\end{conjecture}
It follows quickly from \cite[Corollary 1.5]{green-tao-bszemfrei} and Ruzsa's covering lemma (see e.g. \cite[Lemma 2.14]{tao-vu}) that $G(K) \ll K^{O(\sqrt{K})}$. This is the best upper bound currently known for $G(K)$.

The PFR conjecture allows one to pass from the ``combinatorial information'' $|A + A| \leq K|A|$ to the ``algebraic information'' that $A$ is contained in a union of cosets and back again with only polynomial losses in the constants. By contrast the use of results such as Theorem \ref{mainthm} entail the loss of an exponential, and this is necessary in view of Lemma \ref{ind-pt-example}.  Thus we see that the PFR conjecture would provide significantly more precise information on sets of small doubling constant than Theorem \ref{mainthm}.

We are not able to prove the PFR conjecture in full generality here; however we will be able to verify this conjecture in the special case when $A$ is a \emph{downset}.  
To discuss this model setting we let $(e_i)_{i \in I}$ be the standard basis of $\F_2^I$, and let $\langle, \rangle: \F_2^I \times \F_2^I \to \F_2$ be the standard bilinear form for which $\langle e_i, e_j \rangle = \delta_{ij}$. If $x \in \F_2^I$ then we will write $x_i := \langle x, e_i \rangle \in \F_2$ for the $i$th coordinate of $x$, thus $x = \sum_{i \in I} x_i e_i$. We will write $x \perp e_i$ if $x_i = 0$.

\begin{definition}[Shift-minimal downsets]
Let $A \subseteq \F_2^I$. We say that $A$ is a \emph{downset} if, whenever $i \in I$, $x \perp e_i$ and $x + e_i \in A$, we also have $x \in A$. We say that $A$ is \emph{shift-minimal} if, whenever $i < j$ are elements of $I$, $x \perp e_i$, $x \perp e_j$ and $x + e_j \in A$ we have $x + e_i \in A$. If $A$ is both a downset and shift-minimal then we say that it is a \emph{shift-minimal downset} or \emph{SMD} for short.
\end{definition}

\begin{remarks}
One might usefully think of elements of $\F_2^n$ as $n$-dimensional vectors with coefficients in $\{0,1\}$. Then $A \subseteq \F_2^n$ is a downset if it is closed under replacing ones by zeros, and it is shift-minimal if it is closed under shifting ones to the left.

The terminology comes from the language of set systems, which motivated much of our work in this paper. Given a set $X \subseteq I$ we may associate an element $x_X = \sum_{j \in X} e_j$ of $\F_2^I$. This is clearly a bijection: if $x \in \F_2^I$ we write $X_x = \{ i \in I: x_i = 1 \}$ for the subset of $I$ associated to $x$. In this way we see that there is a one-to-one correspondence between sets $A \subseteq \F_2^I$ and set systems on $I$, that is to say collections of subsets of $I$. We will not make much use of the set system language in this paper, but some readers may find it helpful to think in those terms.
\end{remarks}
 
We can specialise Definition \ref{fg-def} to downsets or SMDs as follows.

\begin{definition}[$\tF(K)$, $\dF(K)$, $\tG(K)$, and $\dG(K)$] Let $K \geq 1$.  
We define $\tF(K)$ to be the least constant such that for any finite $I \subseteq \Z^+$ and any non-empty SMD $A \subseteq \F_2^I$ with doubling constant $\sigma[A] \leq K$, there exists an affine subspace $V \subseteq \F_2^I$ of cardinality $|V|\leq \tF(K) |A|$ which contains $A$.  

We define $\tG(K)$ to be the least constant such that for any finite $I \subseteq \Z^+$ and any non-empty SMD $A \subseteq \F_2^I$ with doubling constant $\sigma[A] \leq K$, there exists a linear subspace $V \subseteq \F_2^I$ of cardinality $|V| \leq |A|$ such that $A$ is covered by at most $\tG(K)$ translates of $V$.

We define $\dF(K)$ and $\dG(K)$ similarly to $\tF(K)$ and $\tG(K)$ but replace ``SMD'' with ``downset'' throughout.
\end{definition}
It is clear that $\tF(K) \leq \dF(K) \leq F(K)$ and that $\tG(K) \leq \dG(K) \leq G(K)$. In fact, in \S \ref{compress-sec} we shall use compression arguments to show
\begin{lemma}[Freiman's theorem is reducible to SMDs]\label{smd}  We have $F(K) = \dF(K) = \tF(K)$ for all $K \geq 1$.
\end{lemma}
In particular, in order to prove Theorem \ref{mainthm} it suffices to work in the SMD model. 

Unfortunately we are not currently able to obtain a similar relationship for $G(K)$, $\dG(K)$, and $\tG(K)$.  Nevertheless, we have the following result.  
\begin{theorem}[PFR for SMDs]\label{mainthm2}
We have $K^{1.46601} \ll \tG(K) \ll K^{30}$ for all $K \geq 1$. 
\end{theorem}
In fact we have the following slightly stronger statement.
\begin{theorem}[PFR for SMDs, again]\label{mainthm2-alt}
If $A \subseteq \F_2^n$ is a SMD, then $A$ is contained in a union of at most $O(\sigma[A]^{30})$ translates of $\F_2^{\lfloor \log_2 |A| \rfloor}$.
\end{theorem}
It is clear that the upper bound in Theorem \ref{mainthm2} follows from Theorem \ref{mainthm2-alt} after performing an order-preserving relabeling to reduce to the ambient group $\F_2^I$ to $\F_2^n$ for some $n$. The lower bound comes by considering a Hamming ball whose radius is given by the solution to a certain optimization problem: the details of this computation are given in \S \ref{polyfrei}.

We also have a similar theorem for downsets.
\begin{theorem}[PFR for downsets]\label{maindown}
We have $K^{1.46601} \ll \dG(K) \ll K^{46}$ for all $K \geq 1$. 
\end{theorem}
Again, we have a slightly stronger statement.
\begin{theorem}[PFR for downsets, again]\label{maindown-alt}
If $A \subseteq \F_2^n$ is a SMD, then $A$ is contained in a union of at most $O(\sigma[A]^{46})$ translates of a coordinate subspace of cardinality at most $|A|$.
\end{theorem}
We shall prove Theorems \ref{mainthm2-alt} and \ref{maindown-alt} (and hence Theorems \ref{mainthm2} and \ref{maindown}) in \S \ref{polyfrei}.  It will be clear from the proofs that the exponents 30 and 46 can be improved, but we do not attempt to optimise them here.

The idea of using methods of extremal set theory to look at Freiman's theorem originated from our study of the work of Bollob\'as and Leader \cite{bollobas-leader}. Our treatment of the basic theory of compression operators in the next section is closely based on their work.

\section{Compressions}\label{compress-sec}

In this section we set out some notation for the machinery of compressions and related structures which we will use throughout the paper, and in particular to establish Lemma \ref{smd} at the end of this section.  Throughout this section $I$ is understood to be a finite set of natural numbers.

\begin{definition}[Lex order]  Suppose that $x = \sum_{i \in I} x_i e_i$ and that $y = \sum_{i \in I} y_i e_i$. Then we define the \emph{lexicographical ordering} (or \emph{lex ordering}) $\prec$ by defining $x \prec y$ if and only if, for the largest coordinate $j$ such that $x_j \neq y_j$, we have $x_j < y_j$.  This is clearly a total ordering.  An \emph{initial segment} in $\F_2^I$ is any set of the form $A := \{ x \in \F_2^I: x \preceq y \}$ for some $y \in \F_2^I$. 
\end{definition}

The following lemma is well known.

\begin{lemma}[Adding initial segments of the lex order]\label{sumset-lex}
Suppose that $A$ and $B$ are initial segments of the lex order on $\F_2^I$. Then so is $A + B$.
\end{lemma}

\proof
 This is easily proved by strong induction on $|A| + |B|$. In fact if $|A| = r$ and $|B| = s$ and if $0 \leq r,s < 2^k$ then the function $\HS(r,s) := |A+B|$ satisfies the relations
\[ \HS(2^k + r,s) = \HS(r,s) + 2^k\]
and \[ \HS(2^k + r, 2^k + s) = 2^{k+1}.\]
(The letters HS stand for Hopf and Stiefel, who studied this function in a context arising in differential topology.) Its r\^ole in the additive combinatorics of $\F_2^I$ was observed by Yuvinsky \cite{yuvinsky}. By induction one may confirm (a result of Plagne \cite{plagne}) that
\begin{equation}\boxeq \HS(r,s) = \min_{j \in \N} 2^j \left( \lceil \frac{r}{2^j}\rceil + \lceil \frac{s}{2^j}\rceil -1 \right).\end{equation}

We now define the notion of a compression with respect to a set of indices.

\begin{definition}[$J$-fibres]
Suppose that $A \subseteq \F_2^I$ and $J \subseteq I$.  Then we may factor $\F_2^I = \F_2^{I \setminus J} \times \F_2^J$ in the obvious manner. For each $x \in \F_2^{I \setminus J}$, define the \emph{$J$-fibre} of $A$ at $x$ to be the set
\[ A_x := \{y \in \F_2^J : (x,y) \in A\} \subseteq \F_2^J.\]
\end{definition}

\begin{definition}[$J$-compressions]
Suppose that $J \subseteq I$ and that $A \subseteq \F_2^I$. Then we define the \emph{$J$-compression} $C_J(A)$ of $A$ to be the set obtained by replacing each fibre $A_x$ by $\overline{A}_x$ for all $x \in \F_2^{I \setminus J}$, the initial segment in the lex order on $\F_2^J$ which has the same cardinality as $A_x$.  We say that $A$ is \emph{$J$-compressed} if $C_J(A)=A$.  If $0 \leq r \leq |I|$, we say that $A$ is \emph{$r$-compressed} if $A$ is $J$-compressed for all $J \subseteq I$ with $|J| \leq r$.  Similarly, we shall refer to a compression operator $C_J$ with $|J| \leq r$ as an \emph{$r$-compression}.
\end{definition}

\begin{example} If $I := \{1,2,3\}$ and $A := \{ e_1, e_1+e_2, e_3\}$, then $C_{\{1\}}(A) = \{ 0, e_2, e_3\}$, $C_{\{2\}}(A) = C_{\emptyset}(A) = A$, $C_{\{1,2\}}(A) = \{0, e_1, e_3\}$, and $C_{\{1,2,3\}}(A) = \{0, e_1, e_2 \}$.
The set $\{0, e_2, e_3\}$ is $1$-compressed but not $2$-compressed.   
\end{example}

The following lemma is trivial:

\begin{lemma}[Trivial properties of compressions]\label{comp-triv}  Let $J \subseteq I$.
\begin{itemize}
\item[(i)] If $A \subseteq \F_2^I$, then $|C_J(A)| = |A|$, and in particular $C_J(A) \subseteq A$ iff $C_J(A) \supseteq A$ iff $C_J(A) = A$.  
\item[(ii)] If $A \subseteq \F_2^I$, then $C_J(C_J(A)) = C_J(A)$.  In other words, $C_J(A)$ is $J$-compressed.
\item[(iii)] If $A \subseteq B \subseteq \F_2^I$, then $C_J(A) \subseteq C_J(B)$.
\end{itemize}
\end{lemma}

When the set $J$ has size at most $2$, or is equal to $I$ one can describe compressions quite explicitly:

\begin{examples}[2-compressions and the full compression]  The compression $C_{\emptyset}$ is trivial: $C_{\emptyset}(A) = A$ for all $A$.  When $J = \{i\}$ the effect of the compression $C_{\{i\}}$ is to move $x + e_i$, where $x \perp e_i$, to $x$ whenever possible. Suppose that $J = \{i,j\}$ with $i < j$ and that $A$ is already a downset. Then the effect of the compression $C_{\{i,j\}}$ is to ``shift'' $x + e_j$, where $x \perp e_i$ and $x \perp e_j$, to $x + e_i$ whenever possible.   If $J=I$, then $C_I(A)$ is simply the initial segment of $\F_2^I$ with length $|A|$.
\end{examples}

From the above description of compressions we immediately obtain

\begin{lemma}[Compressed sets]\label{downsets}
Let $A \subseteq \F_2^I$. 
\begin{itemize}
\item[(i)] $A$ is always $0$-compressed.
\item[(ii)] $A$ is $1$-compressed if and only if it is a downset.
\item[(iii)] $A$ is $2$-compressed if and only if it is an SMD.
\item[(iv)] $A$ is $|I|$-compressed if and only if it is an initial segment of $\F_2^I$.
\end{itemize}
\end{lemma}

Now, we recall the well-known observations that compressions do not increase the size of sumsets.

\begin{lemma}[Compressions and sumsets]\label{compression-sumsets}
Suppose that $A,B \subseteq \F_2^I$ and that $J \subseteq I$. Then $C_J(A) + C_J(B) \subseteq C_J(A + B)$.  In particular $|C_J(A)+C_J(B)| \leq |A+B|$.
\end{lemma}

\proof We follow Bollob\'as and Leader \cite{bollobas-leader}. We proceed by induction on $I$. If $J \neq I$ then we work fibrewise. For any two fibres $A_x, B_y \subseteq \F_2^J$ we have, by the inductive hypothesis, that
\[ \overline{A}_x + \overline{B}_y \subseteq \overline{A_x + B_y}.\]
Thus for any $z \in \F_2^{I \backslash J}$ we have
\begin{align*}
(C_J(A) + C_J(B))_z & = \bigcup_{x + y = z} (\overline{A}_x + \overline{B}_y) \\ & \subseteq \bigcup_{x + y = z} \overline{A_x + B_y} \\ & = \overline{\bigcup_{x + y = z} (A_x + B_y)} \\ & = (\overline{A + B})_z \\ & = C_J(A + B)_z,
\end{align*}
the middle step following from the fact that the initial segments of the lex order are totally ordered under inclusion. It remains to deal with the case $I = J$. By applying the inductive hypothesis and Lemma \ref{comp-term} we see that it suffices to deal with the case in which both $A$ and $B$ are invariant under compressions $C_{J'}$ with $J' \subsetneq I$. 

By an order-preserving relabeling we may take $I = \{1,\ldots,n\}$ for some $n$.
From the compression-invariance assumption we easily confirm that each of $A$ and $B$ is either an initial segment of the lex order or else is equal to the special set 
\[ S := \F_2^{n-1} \setminus \{e_1 + \dots + e_{n-1}\} \cup \{e_n\},\]
a very slight perturbation of an initial segment of lex. Clearly we have $C_{I}(S) = \F_2^{n-1}$. This easily implies that $C_{I}(S) + C_{I}(S) \subseteq C_{I}(S + S)$, so let us assume that $A$ is an initial segment of lex and that $B = S$, so that our task is to prove that
\[ A + \F_2^{n-1} \subseteq C_{I}(A + S).\]
If $A \subseteq \F_2^{n-1}$ then the left-hand side is just $\F_2^{n-1}$ whilst the right contains $C_{I}(S) = \F_2^{n-1}$. If $\F_2^{n-1} \subseteq A$ then the left-hand side is $\F_2^n$. A small amount of thought confirms that $A + S = \F_2^n$ also, thereby concluding the proof.\endproof

As a corollary of Lemma \ref{compression-sumsets} and Lemma \ref{comp-triv}(i) we obtain

\begin{corollary}[Compression decreases doubling constants]\label{compress-double}  Let $J \subseteq I$ and let $A \subseteq \F_2^I$ be non-empty.  Then $\sigma[C_J(A)] \leq \sigma[A]$.
\end{corollary}

Another immediate corollary of Lemma \ref{compression-sumsets} and Lemma \ref{comp-triv}(i) is

\begin{corollary}[Compressed sets closed under addition]\label{downset-add}  Let $J \subseteq I$.  Then the sum of two $J$-compressed sets is still $J$-compressed.  In particular, from Lemma \ref{downsets} we see that the sum of two downsets is a downset, and the sum of two SMDs is an SMD.
\end{corollary}

In fact, for downsets one can express sumsets in terms of lattice operations.

\begin{definition}[Lattice structure] If $x, y \in \F_2^I$ then we write \[ x \vee y := \sum_{i \in I} \max(x_i,y_i) e_i \quad \mbox{and} \quad x \wedge y := \sum_{i \in I} \min(x_i,y_i) e_i.\]  (In the set systems language, these operations correspond to set union and set intersection respectively.)  If  If $A,B \subseteq \F_2^I$ then we write $A \vee B := \{x \vee y : x \in A, y \in B\}$. 
\end{definition}

\begin{lemma}[Lattice representation of downset addition]\label{downlattice}  If $A, B$ are downsets then $A+B= A \vee B$.
\end{lemma}

\proof If $x \in A$ and $y \in B$, we observe that $x+y = x' \vee y'$, where $x' \in A$ and $y' \in B$ are formed from $x$, $y$ by zeroing out those coordinates $i$ for which $x_i=y_i=1$.  This shows that $A+B \subseteq A \vee B$.  The identity $x \vee y = x + y'$ confirms that $A \vee B \subseteq A+B$.
\endproof

The next lemma states that by applying a suitable sequence of compressions to a set $A$ one may obtain a set which is inert with respect to the application of further compressions.

\begin{lemma}[Termination of compression sequences]\label{comp-term}
Let $r$, $1 \leq r \leq |I|$, be an integer and let $A \subseteq \F_2^I$. Then there exists an $r$-compressed set $A'$ which can be obtained from $A$ by applying a finite number of $r$-compressions.
\end{lemma}

\proof By the principle of infinite descent, it suffices to exhibit a weight function $w : \mathcal{P}(\F_2^I) \rightarrow \Z_{\geq 0}$ with the property that $w(C_J(A)) < w(A)$ whenever $C_J$ is an $r$-compression and $C_J(A) \neq A$. It is then clear that one may find a sequence $I_1,\dots,I_k$ with the required properties and with $k \leq w(A)$. There is a wide choice of possible weight functions. One is
\[ w(A) := \sum_{x \in A} 2^{\rho(x)},\]
where $\rho(x)$ is the position of $x$ in the lex order on $\F_2^I$.\endproof

Next, we relate compression to Hamming balls.

\begin{definition}[Hamming balls] If $x \in \F_2^I$ then we define the \emph{Hamming norm} of $x$ as $\Vert x \Vert := |\{ i \in I: x_i = 1 \}|$. If $0 \leq r \leq |I|$ then we define the \emph{Hamming ball} 
\[ H_r(I) := \{x \in \F_2^I : \Vert x \Vert \leq r\}.\]
We abbreviate $H_r(\{1,\ldots,n\})$ as $H_r(n)$.
\end{definition}

\begin{example} We have $H_1(n) = \{0, e_1, \ldots, e_n\}$.
\end{example}

\begin{proposition}[Compression and Hamming balls]\label{comp-prop}  Let $0 \leq r \leq |I|$.
\begin{itemize}
\item[(i)] If $0 \leq r \leq |I|$, then $H_r(I)$ is $2$-compressed.
\item[(ii)] If $H_r(I) \subseteq A \subseteq \F_2^I$ and $C_J$ is a $2$-compression, then $H_r(I) \subseteq C_J(A)$.
\item[(iii)] If $A \subseteq H_r(I)$ and $C_J$ is a $2$-compression, then $C_J(A) \subseteq H_r(I)$.
\end{itemize}
\end{proposition}

\begin{proof} One easily checks that $H_r(I)$ is an SMD, and part (i) then follows from Lemma \ref{downsets}.  Parts (ii) and (iii) then follow from Part (i) and Lemma \ref{comp-triv}(iii).
\end{proof}

\begin{remark}  Observe that the Hamming balls fail to be $3$-compressed in general.  Because of this fact, we will only work with $2$-compressions for the rest of the paper. We have only included the higher compression operators because (for example) Lemma \ref{compression-sumsets} seems to be most naturally proven in this more general context, and because these higher compressions may be of use in future applications.
\end{remark}

Finally, we give some relations between subspaces and the SMD property. 

\begin{lemma}[Compressions and subspaces]\label{smd-subspace}
Suppose that $A \subseteq \F_2^I$ is a SMD and that $H \subseteq \F_2^I$ is an affine subspace of dimension $d$.   Let $I_d$ be the initial segment of $I$ of length $d$.
\begin{itemize}
\item[(i)] For any $J \subseteq I$, $C_J(H)$ is an affine subspace of dimension $d$.
\item[(ii)] $H$ is an SMD if and only if $H = \F_2^{I_d}$.
\item[(iii)] If $H \subseteq A$ then $\F_2^{I_d} \subseteq A$. 
\item[(iv)] If $A \subseteq H$ then $A \subseteq \F_2^{I_d}$.
\item[(v)] If $A = H$ then $A = \F_2^{I_d}$.
\end{itemize}
\end{lemma}
\proof Part (i) can be proven directly. Alternatively, one may use Corollary \ref{compress-double} $C_J(A)$  has doubling constant $1$ and is therefore an affine subspace, which by cardinality considerations must have dimension $d$.

The ``if'' direction of (ii) is trivial.  To prove the ``only if'' direction, observe from rank considerations that $H$ must have a full projection onto $\F_2^J$ for at least one set $J \subseteq I$ with $|J|=d$.  In other words, $H$ contains an element with at least $d$ non-zero entries.  Using the SMD property this implies that $H$ contains $\sum_{i \in I_d} e_i$, and using the downset property we conclude that $H$ contains $\F_2^{I_d}$.  By cardinality or dimension considerations we conclude $H = \F_2^{I_d}$ as desired.

Parts (iii) and (iv) follow from parts (i), (ii) and Lemma \ref{comp-term} (as well as Lemma \ref{downsets}(iii)).  Part
(v) follows immediately from parts (iii) and (iv).
\endproof

We are now in a position to quickly prove Lemma \ref{smd}.\vspace{9pt}

\noindent\emph{Proof of Lemma \ref{smd}.}  The claim $\tF(K) \leq \dF(K) \leq F(K)$ is trivial, so it suffices to show that $F(K) \leq \tF(K)$.  Suppose that $A \subseteq \F_2^I$ has doubling at most $K$; we wish to enclose $A$ in an affine subspace of cardinality at most $\tF(K) |A|$.  By translation we may make $A$ contain the origin, in which case the affine subspace is a linear subspace.  The claim is then equivalent to the assertion that the linear space spanned by $A$ has dimension at most $\log_2( \tF(K) |A|)$.

Suppose that the claim failed.  Then we could find linearly independent elements $v_1,\ldots,v_n$ in $A$ with $n > \log_2( \tF(K) |A|)$, and such that $A$ is contained in the span of $v_1,\ldots,v_n$.  By a linear transformation (and a relabeling of $I$) we may assume that $I \supseteq \{1,\ldots,n\}$ and $v_i = e_i$ for all $i=1,\ldots,n$.  Note that $A$ is contained in the span of $e_1,\ldots,e_n$, i.e. in $\F_2^n$, and so without loss of generality we may take $I = \{1,\ldots,n\}$.  Thus $A$ now contains the Hamming ball $H_1(n)$.  

Applying Lemma \ref{comp-term} (and Lemma \ref{downsets}(iii)) we may form an SMD $A'$ by applying a finite number of $2$-compressions to $A$.  By Lemma \ref{smd-subspace}(i) we see that $A'$ still contains $H_1(n)$.  By Corollary \ref{compress-double} we have $\sigma[A'] \leq \sigma[A] \leq K$.  By Definition \ref{fg-def} we thus see that $A'$ is contained in an affine space $H$ of dimension at most $\log_2( \tF(K) |A| )$.  Applying Lemma \ref{smd-subspace}(iv) we see that $A'$ is contained in a space $\F_2^m$ for some $m \leq \log_2(\tF(K) |A|) < n$.  But this contradicts the fact that $A'$ contains $H_1(n)$.
\endproof

\section{Sumsets and Hamming balls}

Throughout this section we fix a finite set $I$ of natural numbers.
We now use compression technology to obtain some relationships between sumsets and Hamming balls for arbitrary subsets of $\F_2^I$ (not necessarily downsets or SMDs). These results, though somewhat diverse, all confirm that the Hamming ball $H_r(n)$ is a very unstructured set from the additive point of view. We begin with an application of $1$-compressions which shows that large subsets of Hamming balls must expand under sumsets.

\begin{proposition}[Expansion in Hamming balls]\label{hamming-expand}
Suppose that $A,B \subseteq H_r(I)$ for some $r \geq 0$. Then $|A + B| \geq 9^{-r}|A||B|$.
\end{proposition}

\begin{proof}  We are indebted to Oliver Riordan for pointing out the following simple proof; our original proof was significantly more complicated.

Applying $1$-compressions to $A$ and $B$ (using Lemma \ref{compression-sumsets} and Lemma \ref{comp-term}) it suffices to consider the case when $A$ and $B$ are $1$-compressed, and hence (by Lemma \ref{downsets}(ii)) downsets.  By Lemma \ref{downlattice} it now suffices to show that $|A \vee B| \geq 9^{-r}|A||B|$.

To this end, suppose that $x \in A \vee B$. Since $x \in H_{2r}(I)$, the number of representations of $x$ as $a \vee b$ is no more than 
\[ \sum_{i = 0}^{2r} 2^i \binom{2r}{i} = 9^r.\]
The result follows immediately.
\end{proof}

Next, we use $2$-compressions to show that sumsets of large sets must eventually escape Hamming balls.

\begin{proposition}[Expansion out of Hamming balls]\label{downset-bog}
Suppose that $A, A', A'' \subseteq \F_2^I$ is such that 
\begin{equation}\label{aaa}
 |A||A'||A''| \geq (\sqrt{5}+2)^{-r} 2^{3|I|}
\end{equation}
for some $r \geq 0$. Then 
$A + A' + A'' \not \subseteq H_{|I|-r}(I)$.  In other words, $A+A'+A''$ contains a vector with at least $|I| - r$ ones.
\end{proposition}

\proof By an order-preserving labeling we may assume $I = \{1,\ldots,n\}$.
Suppose for contradiction that we can find $A, A', A''$ such that
\begin{equation}\label{eq107} A + A' + A'' \subseteq H_{n-r}(n).\end{equation} 

Applying $2$-compressions to $A, A', A''$ (using Lemma \ref{compression-sumsets}, Lemma \ref{comp-term}, and Lemma \ref{comp-prop}(i)) we see that we can reduce to the case when $A, A', A''$ are all SMDs.  By Lemma \ref{downlattice} (and Corollary \ref{downset-add}), it now suffices to find $x \in A, x' \in A', x'' \in A''$ such that $\Vert x \vee x' \vee x''\Vert > n-r$.

We therefore explore conditions under which there are $x \in A, x' \in A', x'' \in A''$ for which $\Vert x \vee x' \vee x''\Vert$ is large. Let us imagine that $\Vert x\Vert + \Vert x' \Vert + \Vert x'' \Vert$ counters have been placed on the points of $I$, the number of counters placed on $i$ being $x_i + x'_i + x''_i$. Let us begin at $i = n$ and try to redistribute these counters by moving them to the left in such a way that as many elements of $I$ are covered by at least one counter. Thus for $i = n, n-1, n-2,\dots$ in turn we pick up all of the counters at that place, and then put one counter down if we have one available in our hand. Write $S$ for that set of points which have at least one counter at the end, and set $z := \sum_{i \in S} e_i$. By construction we see that
\[ z = \overline{x} \vee \overline{x'} \vee \overline{x''},\]
where $\overline{x}$ is obtained from $x$ by removing some ones from the coordinate vector of $x$ and shifting others to the left, and similarly for $\overline{x'}$ and $\overline{x''}$.
Since $A, A', A''$ are SMDs it follows that $\overline{x} \in A$, $\overline{x'} \in A'$ and $\overline{x''} \in A''$, and hence that $z \in A + A' + A''$.

We define the \emph{defect} $d(x,x',x'')$ to be the number of positions which have no counters when our process is finished, that is to say $d(x,x',x'') := n - |S| = n - \Vert z \Vert$.

Our aim is to estimate the number of triples $(x,x',x'') \in (\F_2^n)^3$ for which $d(x,x',x'') \geq r$. To do this, we first imagine a probabilistic variant of the above ``moving counters to the left'' game which takes place on $\Z_{\leq 0}$ instead of $[n]$. The number of counters on each point $0, -1, -2, \ldots$ is chosen by tossing three coins, and placing a counter for each coin that comes up heads. All coin tosses are independent. Once more one imagines moving the counters to the left so as to cover as many points as possible. For $n \geq 0$ we define the $n$-\emph{defect} $D_n$ (which is a random variable) to be the number of positions among $\{-n,-(n-1),\dots,-1,0\}$ which are left uncovered by our process of moving counters to the left. Define $D = \sup_{n \geq 0}D_n$; we will see shortly that this quantity is almost surely finite. Set $q(r) := \P(D \geq r)$, thus for instance $q(0)=1$.

Suppose we play the same game but starting from $n=-1$ instead of $n=0$, ignoring the $0$ position entirely.  The analogous random variable $D'$ has the same probability distribution as $D$, but obeys the relationship
$$ D \geq r+1 \hbox{ if and only if } D' \geq r+i_0$$
for any fixed $r \geq 0$, where $i_0$ is the number of counters at the position $0$.  Conditioning this relationship on the random variable $i_0$, we obtain the recurrence relation
\begin{equation}\label{recurrence} q(r+1) = \textstyle\frac{1}{8} q(r+3) + \frac{3}{8}q(r+2) + \frac{3}{8}q(r+1) + \frac{1}{8}q(r)\displaystyle \end{equation}
for all $r \geq 0$. 

By conditioning upon the position of the $r^{\operatorname{th}}$ uncovered counter after the completion of our process of moving counters to the left, one may also confirm the multiplicativity property
\[ q(r+s) = q(r)q(s)\]
for all $r,s \geq 1$ (note that the claim is trivial if $r=0$ or $s=0$).  Thus $q(r) = q(1)^r$.
The characteristic polynomial for the linear recurrence \eqref{recurrence}, $x^3 + 3x^2 - 5x + 1$, factors as $(x-1)(x^2 + 4x - 1)$ and so (since one clearly has $q(1) < 1$) it follows easily that 
\begin{equation}\label{q-formula} q(r) = (2 + \sqrt{5})^{-r}.\end{equation}
Now if $x,x',x''$ are chosen independently at random from $(\F_2^n)^3$ then 
\[ \P(d(x,x',x'') \geq r) \leq q(r) = (2 + \sqrt{5})^{-r}.\]
The result follows.\endproof

\begin{remarks} Our treatment is essentially taken from \cite{frankl-survey}. In that article upper bounds are obtained for the maximum size of $3$-wise $t$-intersecting set systems: set systems $\mathcal{F},\mathcal{F}',\mathcal{F}''$ such that $|F \cap F' \cap F''| \geq t$ for all $F \in \mathcal{F}$, $F' \in \mathcal{F}'$, $F'' \in \mathcal{F}''$. Equivalently we have $|F^c \cup F^{\prime c} \cup F^{\prime\prime c}| \leq n - t$. This condition is very similar to \eqref{eq107}, and is in fact \emph{equivalent} to it in the case that $\mathcal{F},\mathcal{F}'$ and $\mathcal{F}''$ are downsets.
\end{remarks}

An equivalent formulation of Proposition \ref{downset-bog} is that if we have the largeness condition \eqref{aaa}, then $A + A' + A'' + H_r(I) = \F_2^I$.  We can now ``bootstrap'' this proposition to the case in which $A, A', A''$ do not have to obey \eqref{aaa}.

\begin{proposition}[Expansion result]\label{expansion-prop}
Suppose that $A,A',A'' \subseteq \F_2^I$.
Then
\[ |A + A' + A'' + H_r(I)| \geq \min\left(\half(2 + \sqrt{5})^{r/3} |A|^{1/3}|A'|^{1/3}|A''|^{1/3}, 2^{|I|}\right).\]
\end{proposition}
\proof By relabeling we may take $I = \{1,\ldots,n\}$.
Write
\[ f(n) := \min_{\substack{A,A',A'' \subseteq \F_2^n \\ |A||A'||A''| \leq 8(2+\sqrt{5})^{-r}2^{3n}}} \frac{|A + A' + A'' + H_r(n)|}{|A|^{1/3}|A'|^{1/3} |A''|^{1/3}}.\]

\begin{claim} For all $n$ we have
\begin{equation}\label{ind-hyp} f(n) \geq \half(2 + \sqrt{5})^{r/3}.\end{equation}
\end{claim}

\begin{proof}  We induct on $n$.  Applying $2$-compressions together with Lemma \ref{comp-term}, Lemma \ref{compression-sumsets} and Proposition \ref{comp-prop}(i), we see that the minimum is attained when $A,A',A''$ are 2-compressed, and thus (by Lemma \ref{downsets}(iii)) SMDs. It may be that the minimum is attained for sets $A,A',A''$ with $|A||A'||A''| \leq (2 + \sqrt{5})^{-r} 2^{3n}$. Otherwise (of course) the minimum is attained for sets with $|A||A'||A''| \geq (2 + \sqrt{5})^{-r}2^{3n}$. We call these cases 1 and 2 respectively (these cases can overlap if the minimum is achieved multiple times).

Suppose that we are in case 1.
Split $\F_2^n$ as $\F_2^{n-1} \times \F_2$ in some arbitrary direction and write $A_i := \{x \in \F_2^{n-1}: (x,i) \in A\}$ for $i = 0,1$. Define $A'_i, A''_i$ similarly. Set $p := |A_0|/|A|, p' := |A'_0|/|A|$ and $p'' := |A''_0|/|A|$. The downset property implies that $p,p',p'' \geq 1/2$; assume without loss of generality that $p'' \geq p' \geq  p \geq 1/2$. 

Observe that $A + A' + A'' + H_r(n)$ contains disjoint copies of $A_0 + A'_0 + A''_0 + H_r(n-1)$ and $A_1 + A'_0 + A''_0 + H_r(n-1)$. Note also that
\[(p'p'')^{1/3} \left( p^{1/3} + (1 - p)^{1/3} \right)  \geq p + p^{2/3}(1 - p)^{1/3}  \geq  1.
\] Since we have
\[ |A_{i}||A'_{i'}||A''_{i''}| \leq 8(2 + \sqrt{5})^{-r} 2^{3(n-1)}\] for any choice of $i,i',i'' \in \{0,1\}$, we see that $f(n) \geq f(n-1)$ in this case, and the claim follows by induction hypothesis.

In case 2 we are in the situation covered by Proposition \ref{downset-bog}. It follows from that proposition that $A + A' + A''$ contains a vector with at least $n - r$ ones. It follows that $A + A' + A'' + H_r(n)$ contains the vector $e_1 + \dots + e_n$ and hence, since $A + A' + A'' + H_r(n)$ is an SMD, that
\[ A + A' + A'' + H_r(n) = \F_2^n.\] Thus
\[ f(n) \geq \frac{2^n}{|A|^{1/3}|A'|^{1/3}|A''|^{1/3}} \geq \half(2 + \sqrt{5})^{r/3},\] as claimed in \eqref{ind-hyp}.
\end{proof}

We return now to the proof of Proposition \ref{expansion-prop}.  If $|A||A'||A''| \leq 8(2+\sqrt{5})^{-r} 2^{3n}$, then we have
\[ |A + A' + A'' + H_r(n)| \geq \half(2 + \sqrt{5})^{r/3}|A|^{1/3}|A'|^{1/3}|A''|^{1/3}.\]
If instead we have $|A||A'||A''| < 8(2+\sqrt{5})^{-r} 2^{3n}$ we may apply Proposition \ref{downset-bog} to conclude that
\[ A + A' + A'' + H_r(n) = \F_2^n.\]
In either case we conclude Proposition \ref{expansion-prop}.\endproof

\section{The polynomial Freiman-Ruzsa conjecture for SMD sets and downsets}\label{polyfrei}

We use the above machinery, together with the standard tools of Pl\"unnecke-Ruzsa theory (see e.g. \cite[Chapter 2]{tao-vu}) to prove Theorem \ref{mainthm2}, or equivalently Theorem \ref{mainthm2-alt}.  We first prove a slight variant of this theorem, which will be useful for establishing all of our main theorems.

\begin{proposition}[PFR for SMD sets, intersection version]\label{prop7.1}
Suppose that $A \subseteq \F_2^m$ is a SMD and that $\sigma[A] \leq K$. Set $m := \lfloor \log_2|A|\rfloor$. Then 
\[ |A \cap \F_2^{m}| \gg K^{-29}|A|.\]
\end{proposition}
\proof
Set $m' := m + 5\log_2 K + 1$.  By enlarging $n$ if necessary we may assume that $n = m'+t$ for some $t \geq 0$.  Factor $\F_2^{m'+t} = \F_2^{m'} \times \F_2^{\{m'+1,\dots,m'+t\}}$ and consider, for each $x \in \F_2^{m'}$, the fibre $A_x$ of $A$ above $x$. Write $B = A \cap \F_2^{m'}$ for the ``base'', and let $\pi : \F_2^{m' + t} \rightarrow \F_2^{m'}$ be the projection map.  

Suppose that $A_0$ contains a vector with Hamming norm $r$.
Since $A$ is a SMD, we conclude that $H_r(m') \subseteq B$, and therefore by Proposition \ref{expansion-prop} we have
\begin{equation}\label{bubba}
|B + B + B + B| \geq |B + B + B + H_r(m')| \geq \min\left(\half (2 + \sqrt{5})^{r/3}|B|, 2^{m'}\right).
\end{equation}
On the other hand, since $A$ is a downset, we have $A_x \subseteq A_0$ for all $x \in B$, and so
$$ |A| \leq |A_0| |B|.$$
Also, observe that $B+B+B+B+A_0 \subseteq A+A+A+A+A$, and so
$$ |B+B+B+B| |A_0| \leq |A+A+A+A+A|.$$
On the other hand, from the Pl\"unnecke-Ruzsa inequality (see \cite[Corollary 6.29]{tao-vu}) and the hypothesis $|A+A| \leq K|A|$ we have
$$ |A+A+A+A+A| \leq K^5 |A|$$
From the preceding three inequalities\footnote{The observant reader will see that this argument can in fact be easily modified to the case when $A$ is not a downset, simply by replacing $A_0$ with $A_x$ where $x$ maximises $|A_x|$.} we conclude
$$ |B+B+B+B| \leq K^5 |B|.$$
The choice of $m'$ guarantees that $K^5 |B| \leq K^5|A| < 2^{m'}$, and so we are forced to conclude from \eqref{bubba} that
\[ \half (2 + \sqrt{5})^{r/3} \leq K^5\] 
and hence that
$$ r \leq 15 \log_{2 + \sqrt{5}} K + 2.$$
In other words, $A_0$ is contained in $H_r(\F_2^{\{m+1,\dots,m+t\}})$, where $r := \lceil 15 \log_{2+\sqrt{5}} K + 2\rceil$. Applying Proposition \ref{hamming-expand} we have 
\[ |(A + A)_x| \geq |A_0 + A_x| \geq 9^{-r}|A_0||A_x|,\] and therefore, summing over $x$, we see that
\[ |A+A| \geq 9^{-r}|A_0||A|.\] It follows immediately that
\[ |A_0| \leq 9^rK \ll K^{1 + 15\log_{2+\sqrt{5}} 9} \leq K^{24}.\]
Now $A$ is a downset, and so $A_x \subseteq A_0$ for every $x$. It follows immediately that
\[ |B| \gg K^{-24}|A|.\] 
Recall that $B := A \cap \F_2^{m'}$. The fact that $A$ is a downset easily implies that 
$|A \cap \F_2^j|/2^j$ is a non-increasing function of $j \in \N$. Recalling the definition of $m'$ in terms of $m$ it follows that
\[ |A \cap \F_2^m| \gg K^{-29}|A|,\] as claimed.  This proves Proposition \ref{prop7.1}, and hence the upper bound of Theorem \ref{mainthm2}.
\endproof

\noindent\emph{Proof of Theorem \ref{mainthm2-alt}.}  Set $K := \sigma[A]$, $m := \lfloor \log_2 |A| \rfloor$, and $B := A \cap \F_2^m$.  Applying Ruzsa's covering lemma (see e.g. \cite[Lemma 2.14]{tao-vu}) we may cover $A$ by at most $\frac{|A+B|}{|B|}$ translates of $B-B$.  But $|A+B| \leq |A+A| \leq K|A|$ and $B-B \subseteq \F_2^m$, while from Proposition \ref{prop7.1} we have
$|B| \gg K^{-29} |A|$.  The claim follows.
\endproof

Now we turn from SMD sets to downsets.  The key proposition is

\begin{proposition}[Small doubling implies escape from Hamming ball]\label{pprop}  Let $A \subseteq \F_2^n$ be such that with $\sigma[A] \leq K$, and such that $A+A+A \subseteq H_r(n)$.  Then $r \geq \log_2 |A| - 42 \log_2 K - O(1)$.
\end{proposition}

\proof By applying $2$-compressions (using Lemma \ref{compression-sumsets}, Corollary \ref{compress-double}, Lemma \ref{comp-term} and Proposition \ref{comp-prop}) we may reduce to the case where $A$ is $2$-compressed, and hence (by Lemma \ref{downsets}(iii)) an SMD set.  

Let $m := \lfloor \log_2 |A| \rfloor$, and set $B := A \cap \F_2^m$.  Applying Proposition \ref{prop7.1} we see that
$$ |B| \gg K^{-29}|A| \gg K^{-29} 2^m.$$
Applying Proposition \ref{downset-bog} we see that $B+B+B$ contains a vector of norm at least
$$ m - 3 \log_{\sqrt{5}+2} K^{29} - O(1) \geq \log_2 |A| - 42 \log_2 K - O(1).$$
But $B+B+B$ is contained in $A+A+A$, which lies in $H_r(n)$ by hypothesis.  The claim follows.
\endproof

\noindent\emph{Proof of Theorem \ref{maindown-alt}.}  Let $K := \sigma[A]$.  Let $A$ be a downset with $\sigma[A] \leq K$.  By Proposition \ref{pprop} $A+A+A$ contains a vector with Hamming norm at least $\log_2 |A| - 42 \log_2 K - O(1)$.  By Corollary \ref{downset-add}, $A+A+A$ is a downset, and so $A+A+A$ contains a coordinate subspace $\F_2^J$ of cardinality $\gg K^{-42} |A|$. Now observe that
$$ |A+\F_2^J| \leq |A+A+A+A| \leq K^4 |A| \ll K^{46} |\F_2^J|$$
by the Pl\"unnecke-Ruzsa inequality (see \cite[Corollary 6.29]{tao-vu}).  
Thus $A$ can be covered by $O(K^{46})$ translates of $\F_2^J$.  If $|\F_2^J| \leq |A|$ we are done, so suppose instead that $|\F_2^J| > |A|$.  Then we have
$$ |\F_2^J| \leq |A+\F_2^J| \leq |A+A+A+A| \leq K^4 |A| \leq K^4 |\F_2^J|$$
so that $A$ can be covered by $O(K^4)$ translates of $\F_2^J$, which can in turn be covered by $O(K^4)$ translates of a vector space of cardinality at most $|A|$.  Thus in both cases we obtain Theorem \ref{maindown-alt}, and the upper bound of Theorem \ref{maindown} follows.
\endproof

It remains to establish the lower bound in Theorems \ref{mainthm2} and \ref{maindown}. Since $\tG(K) \leq \dG(K) \leq G(K)$, it is clear that the lower bound in Theorem \ref{mainthm2} is stronger.\vspace{11pt}

\noindent\emph{Proof that $\tG(K) \gg K^{1.46601}$.}
Let $\alpha \in (0,1/10]$ be a constant to be optimised later. Let $A = H_{\lfloor \alpha n\rfloor}(n) \subseteq \F_2^n$; this is of course a SMD. A standard computation using Stirling's formula shows that
\[ |A| = 2^{(h(\alpha) + o(1))n},\]
where $h$ is the entropy function $h(\alpha) := -\alpha \log_2\alpha - (1-\alpha)\log_2(1-\alpha)$. Since $A + A = H_{2\lfloor \alpha n\rfloor}(n)$ we have
\begin{equation}\label{eq331} \frac{1}{n}\log_2 \sigma[A] = h(2\alpha) - h(\alpha) + o(1).\end{equation}
Set $m := \lfloor \log_2|A|\rfloor$. Since $A$ is a SMD, we see from Lemma \ref{smd-subspace} that if $A$ is covered by $t$ translates of a subspace of dimension $m$ then it is in fact covered by $t$ translates of $\F_2^m$. Thus if $\tau[A]$ is the minimum possible value of $t$ then we see that
\[ \tau[A] = |H_{\alpha n}(n - m)|\]
and hence that
\begin{equation}\label{eq332} \frac{1}{n}\log_2 \tau[A] = (1 - h(\alpha))h\big( \frac{\alpha}{1-h(\alpha)}\big) + o(1).\end{equation}
Using \emph{Mathematica} we found that the ratio between \eqref{eq332} and \eqref{eq331} is largest when $\alpha = 0.0939288\ldots$, in which case it is $1.46601\ldots + o(1)$.\endproof

\begin{remarks}
Proposition \ref{downset-bog} immediately implies that if $A \subseteq \F_2^n$ is a downset and if $|A| = \alpha 2^n$ then $A + A + A$ contains a subspace of dimension at least $n - O(\log(1/\alpha))$. Proposition \ref{pprop} implies a similar conclusion under the assumption that $A$ is a downset with $\sigma[A] \leq K$: now $A + A + A$ contains a subspace of dimension at least $\log_2|A| - O(\log K)$. Results of a similar type, but with weaker dependencies, are known without the assumption that $A$ is a downset. If now $A$ is an arbitrary set and $|A| = \alpha 2^n$ then $A + A + A$ contains an affine  subspace of dimension $n - O(\alpha^C)$, whilst if $\sigma[A] \leq K$ then $A + A + A$ contains an affine subspace of dimension $\log_2|A| - O(K^C)$. The first of these results is essentially due to Bogolyubov \cite{bog} (see also \cite[Proposition 4.39]{tao-vu}), whilst the second follows from the first together with Ruzsa's observation \cite{ruzsa-freiman-old} that a set with doubling $\sigma[A] \leq K$ is Freiman $6$-isomorphic to a subset $A' \subseteq \F_2^m$, where $|A'| \geq K^{-12}2^m$ (cf. \cite[Lemma 5.26]{tao-vu}, as well as \cite{green-finite-field-survey}).

It might reasonably be conjectured that the bounds that we have obtained in the downset case hold in general. Thus it is possible that if $A \subseteq \F_2^n$ has size $|A| = \alpha 2^n$ then $A + A + A$ contains an affine subspace of dimension at least $n - O(\log(1/\alpha))$. This conjecture, which we call the \emph{Polynomial Bogolyubov Conjecture}, implies the PFR Conjecture, but seems to be strictly stronger than it. Our results verify this conjecture in the downset model.
\end{remarks}

\section{Freiman-type theorems in $\F_2^n$}\label{freisec}

We now turn our attention to the proof of Theorem \ref{mainthm}.  We first need a preliminary result:

\begin{proposition}[Freiman's lemma in finite fields]\label{freiman-lemma-cube}
Suppose that $A, B \subseteq \F_2^n$ satisfy $H_1(n) \subseteq B \subseteq A$. Suppose also that $|A| \leq n^C$ for some $C \geq 1$. Then 
\[ |A + B| \geq \left(\half n - O(\sqrt{n} \log n)\right)|B|,\] where the implied constant depends only on $C$.
\end{proposition}

\begin{remark}
Similar results, with the $O(\sqrt{n \log n})$ replaced by something weaker such as $o(n)$, are available under much weaker hypotheses, for example that $|A| \leq 2^{n^{1-\delta}}$ for some $\delta > 0$.
This proposition should be compared with
Freiman's lemma \cite{freiman-lemma}, which asserts that if $A$ is a subset of $\R^d$ which is not contained in any proper affine subspace then 
\[ |A+A| \geq (d+1)|A| - \half d(d+1).\]
\end{remark}

\proof We may assume that $n$ is large, since the claim is vacuously true otherwise.
Suppose, as a hypothesis for contradiction, that the result is false. By applying $2$-compressions (using Lemma \ref{comp-term}, Lemma \ref{comp-triv}(iii), Lemma \ref{compression-sumsets}, and Proposition \ref{comp-prop}) we may assume that $A$ and $B$ are SMDs. The fact that $A$ is a downset, together with the upper bound on $|A|$, implies that every element of $A$ has Hamming norm at most $C \log_2 n$, or in other words $A \subseteq H_{C\log_2 n}(n)$. Set $M := \lfloor 10C \sqrt{n} \log_2 n \rfloor$. 

\begin{claim} If $x \in A$ there cannot be more than one value of $i \geq M$ such that $x_i = 1$.
\end{claim}

\begin{proof}
 Suppose that there is some $x \in A$ such that there are two values of $i \geq M$ such that $x_i = 1$. By shifting to the left, we see that $H_2(M) \subseteq A$. We bound from below the number of distinct elements $z = a+b$ with $a \in H_2(M)$, $b \in B$ and $a \wedge b = 0$. 

On the one hand the number of choices for the pair $(a,b)$ is at least
\[ |B| \binom{M - C\log_2 n}{2} \geq (8 + o(1))C^2|B| n \log^2_2 n.\]

On the other hand, no $z$ arises in too many ways as such a sum $a + b$. Indeed
since $\Vert z \Vert \leq 2C\log_2 n$, the number of choices for $a$ is no more than $(2C\log_2 n)^2$. 

Comparing these two bounds it follows that
\[ |A+B| \geq (8 + o(1))C^2(2C\log_2 n)^{-2} |B| n\log^2_2 n > (1 - o(1)) n |B|.\] 
This is contrary to assumption (for $n$ large enough), and so the claim follows.
\end{proof}

For each $i$, $M \leq i \leq n$, write $B_i := \{x \in B: x_i = 1\}$. Set $B_0 := B \setminus \bigcup_i B_i$. Now $A + B$ contains the sets
\[ e_i + B_0 : M \leq  i\leq n\]
and
\[ e_i + B_j : M \leq i,j \leq n, i \neq j.\]
No element lies in more than two of these sets, and therefore
\[ |A+B| \geq \frac{1}{2}\left( (n-M)|B_0| + (n - M - 1)\sum_{j = M}^n |B_j|\right) \geq \frac{n - M - 1}{2}|B|.\]
This concludes the proof of Proposition \ref{freiman-lemma-cube}.\endproof

\noindent\emph{Proof of Theorem \ref{mainthm}.}
In view of Lemma \ref{ind-pt-example} and Lemma \ref{smd}, it suffices to show that $\tF(K) \leq 2^{2K + O(\sqrt{K} \log K)}$.  

Let $A$ be an SMD, and let $n$ be the largest integer such that $H_1(n) \subseteq A$.  Since $A$ is a downset, we see that $A \subseteq \F_2^n$.  It will thus suffice to show that $n \leq \log_2 |A| + 2K + O(\sqrt{K} \log K)$.  

Let $m := \lfloor \log_2 |A| \rfloor$ and factor $\F_2^n = \F_2^m \times \F_2^{n'}$ where $n' := n - m$. For each $x \in \F_2^{m}$ write $A_x$ for the fibre 
\[ A_x := \{ y \in \F_2^{n'}: (x,y) \in A\}.\]
Let $B$ be the ``base'':
\[ B := A \cap \F_2^m = \{ x: A_x \neq \emptyset\}.\]
Since $A$ is a downset, we have $A_x \subseteq A_0$ for all $x$. Now $A + A$ contains $B + A_0$, and therefore
\[ |B||A_0| = |B + A_0| \leq K|A|.\]
On the other hand from Proposition \ref{prop7.1} we have
\[ |B| \gg K^{-29}|A|.\]
Comparing these two bounds we see that $|A_0| \ll K^{30}$. If $n' \leq K$ then we are done, and if not the hypotheses of Lemma \ref{freiman-lemma-cube} are satisfied and so 
\[ |A_x + A_0| \geq \left(\frac{n'}{2} - O(\sqrt{n'\log n'})\right) |A_x|\] for all $x$.
It follows immediately that
\begin{align*} K|A| &\geq |A + A|\\
& \geq \sum_{x}|A_x + A_0|\\ & \geq \left(\frac{n'}{2} - O(\sqrt{n'\log n'})\right)\sum_x |A_x|\\ & \geq \left(\frac{n'}{2} - O(\sqrt{n'}\log n')\right)|A|,\end{align*}
and therefore
\[ n' \leq 2K + O(\sqrt{K}\log K).\]
This concludes the proof.
\endproof

\section{The small $K$ case}\label{smallk}

Most of our paper has been concerned with the behaviour of quantities such as $F(K)$ and $G(K)$ in the asymptotic limit $K \to \infty$.  For sake of completeness we make some remarks about the case of small $K$, which involves rather different techniques to those employed above.  

In \cite{desh} it was shown that $F(K) \leq \frac{2K-1}{3K - K^2 - 1}$ for $1 \leq K \leq \frac{12}{5}$, together with a more complicated upper bound for $F(K)$ in the region $\frac{12}{5} < K < 4$.  Here we give a precise formula for $F(K)$ in the region $1 \leq K < \frac{9}{5}$.

\begin{proposition}[$F(K)$ for small $K$]\label{fk-small}  If $K < \frac{7}{4}$, then $F(K) = K$, while $F(K)=\frac{8}{7}K$ for $\frac{7}{4} \leq K < \frac{9}{5}$.  \textup{(}In particular, $F$ is discontinuous.\textup{)}
\end{proposition}

\proof The lower bound $F(K) \geq K$ (for any $K$) is clear by setting $A$ to be a random subset of a large affine space $V$ of density arbitrarily close to $1/K$.  The lower bound $F(K) \geq \frac{8}{7} K$ for $K \geq \frac{7}{4}$ follows by setting $A$ to be a random subset of $\{0, e_1, e_2, e_3\} \times \F_2^{n-3}$ with $n$ large of density arbitrarily close to $\frac{7}{4K}$, where $e_1,\ldots,e_n$ is the standard basis of $\F_2^n$.

Now we show $F(K) \leq K$ for $K < \frac{7}{4}$, using an argument from \cite[Corollary 5.6]{tao-vu}.  
From Kneser's theorem (see e.g. \cite[Theorem 5.5]{tao-vu}) we have $|A+A| \geq 2|A|-|H|$, where $H$ is the largest subspace of $\F_2^n$ with the property that $A+A$ is a union of cosets of $H$.  Since $|A+A| < \frac{7}{4} |A|$, we conclude that $|A+A| < 7|H|$, thus $A+A$ is the union of at most six cosets of $H$.  In particular, if $B := (A+H)/H \subseteq \F_2^n/H$, then $B+B$ has cardinality at most $6$, and $B+B$ cannot be expressed as the union of cosets of any non-trivial subspace of $\F_2^n/H$.  Kneser's theorem again gives $|B+B| \ge 2|B|-1$ and thus $|B| \leq 3$.  If $|B|=2$ or $|B|=3$ then $B+B$ is a subspace of dimension $1$ or $2$ respectively, which is not permitted, so $B$ is a singleton.  Thus $A$ is contained in a coset $V$ of $H$.  Since $|V| = |A+A| \leq K|A|$, we obtain $F(K) \leq K$ as desired.

Now take $\frac{7}{4} \leq K < \frac{9}{5}$.  Repeating the above argments with the same notation we obtain that $|B+B| < 9$ and hence $|B| \leq 4$.  The only new case is when $|B|=4$, in which case $B+B$ is contained in a subspace of dimension $3$ and has cardinality exactly $7$.  This implies that $A$ is contained in an affine subspace of cardinality $8|H|$.  Since $7|H| = |A+A| \leq K |A|$, the bound $F(K) \leq \frac{8}{7} K$ follows.
\endproof

\begin{remark} A similar analysis can (in principle) be used to compute $F(K)$ exactly for any given $K < 2$, though the number of cases to consider via this method grows quickly to infinity as $K \to 2$.
\end{remark}

The same analysis also gives precise values for $G(K)$ when $G$ is small.

\begin{proposition}[$G(K)$ for small $K$]  If $1 < K < \frac{7}{4}$, then $\tG(K) = \dG(K) = G(K) = 2$, while $\tG(K) = \dG(K) = G(K) = 3$ for $\frac{7}{4} \leq K < \frac{9}{5}$.  
\end{proposition}

\begin{proof}  The lower bound $\tG(K) \geq 2$ for $1 < K < \frac{7}{4}$ can be established by setting $A$ to be an initial segment in a large space $\F_2^n$ of density arbitrarily close to $1/K$.  The lower bound $\tG(K) \geq 3$ for $\frac{7}{4} \leq K < \frac{9}{5}$ can be established by setting $A$ to be the Cartesian product (in $\F_2^n \times \F_2^3 = \F_2^{n+3}$) of an initial segment in a large space $\F_2^n$ of density arbitrarily close to $\frac{7}{4K}$, together with the set $\{0, e_{n+1}, e_{n+2}, e_{n+3} \}$.

Since $\tG(K) \leq \dG(K) \leq G(K)$, it thus suffices to establish the bounds $G(K) \leq 2$ for $K < \frac{7}{4}$ and $G(K) \leq 3$ for $K < \frac{9}{5}$.  If $K < \frac{7}{4}$, then by the analysis used to prove Lemma \ref{fk-small} we see that $A$ is contained in an affine space $V$ of size $|V| \leq K|A|$.  Subdividing $V$ into two affine spaces of size at most $K|A|/2 \leq |A|$ we obtain the bound $G(K) \leq 2$.  (The case when $|V|=1$ can of course be dealt with by hand.)

If $K < \frac{9}{5}$, the only new case is when $A$ is the union of four cosets of $H$, where $7|H| \leq K|A|$.  In particular, we see that $2|H| \leq \frac{2}{7} \frac{9}{5} |A| < |A|$. Covering these four cosets of $H$ by three cosets of a subspace of one higher dimension than $H$ we obtain the bound $G(K) \leq 3$.
\end{proof}

\begin{remark} While these above small examples show that $\tG(K) = \dG(K) = G(K)$ for small $K$, we do not expect this relationship to hold exactly in general; it may be that there are non-downsets of reasonably small doubling constant which are significantly harder to cover efficiently by subspaces than downsets or SMDs with the same doubling constant. Certainly the arguments used to prove Lemma \ref{smd} do not seem to extend to the $G()$ family of quantities. On the other hand, Conjecture \ref{pfr} would clearly imply that $\tG(K)$, $\dG(K)$, and $G(K)$ are equivalent up to polynomial factors.
\end{remark}

\end{document}